\documentclass[11pt]{article}
\usepackage{amsmath} \usepackage{amssymb} \usepackage{latexsym} \usepackage{enumitem} \usepackage{mathrsfs} \usepackage{comment}
\usepackage{color}

\usepackage[colorlinks=true,urlcolor=blue,
citecolor=red,linkcolor=blue,linktocpage,pdfpagelabels,bookmarksnumbered,bookmarksopen]{hyperref}

\newtheorem{assumption}{Assumption}

\def\qed{ \ \vrule width.2cm height.2cm depth0cm\smallskip}
\newenvironment{proof}{\noindent {\bf Proof.\/}}{$\qed$\vskip 0.1in}

\newcommand{\ba}{\begin{array}}
\newcommand{\ea}{\end{array}}
\newcommand{\be}{\begin{equation}}
\newcommand{\ee}{\end{equation}}
\newcommand{\bea}{\begin{eqnarray}}
\newcommand{\eea}{\end{eqnarray}}
\newcommand{\beaa}{\begin{eqnarray*}}
\newcommand{\eeaa}{\end{eqnarray*}}

\def\dbR{\mathbb{R}}

\def\a{\alpha}

\def\d{\delta}

\def\cA{{\cal A}}

\def\cX{{\cal X}}

\def\no{\noindent}

\def\ms{\medskip}
\def\bs{\bigskip}
\def\q{\quad}
\def\qq{\qquad}

\def\cds{\cdots}

\def\qed{ \hfill \vrule width.25cm height.25cm depth0cm\smallskip}

\newcommand{\basa}{\begin{assumption}}
\newcommand{\easa}{\end{assumption}}

\newcommand{\bas}{\begin{assum}}
\newcommand{\eas}{\end{assum}}

\def\limsup{\mathop{\overline{\rm lim}}}
\def\liminf{\mathop{\underline{\rm lim}}}

\def\cds{\cdots}

\def\dis{\displaystyle}

\def\1{{\bf 1}}

\def\:{\!:\!}
\def\reff{\eqref}
\def \proof{{\noindent \bf Proof.\quad}}

at 9pt

\definecolor{alp}{rgb}{0.0, 0.5, 0.0}

\newtheorem{thm}{Theorem}[section]

\newtheorem{prop}[thm]{Proposition}
\newtheorem{rem}[thm]{Remark}

\newtheorem{defn}[thm]{Definition}
\newtheorem{assum}[thm]{Assumption}

\begin{document}

\title{\bf Is a sophisticated agent always a wise one?} 
\author{Jianfeng Zhang\thanks{\noindent  Department of Mathematics, 
University of Southern California, Los Angeles, CA 90089, USA. E-mail:
\href{mailto:jianfenz@usc.edu}{jianfenz@usc.edu}. The author is supported in part by NSF grants DMS-1908665 and DMS-2205972.
}  \thanks{\noindent The author would like to thank Jiongmin Yong for very helpful discussions and two anonymous referees for their inspiring comments.}
}
%\date{}
\maketitle

\begin{abstract} For time inconsistent optimal control problems, a quite popular approach is the equilibrium approach, taken by the sophisticated agents. In this short note we construct a deterministic continuous time example where the unique equilibrium is dominated by another control. Therefore, in this situation it may not be wise to take the equilibrium strategy.
\end{abstract}

\no{\bf Keywords.}  Time inconsistency, Pareto optimal, sophisticated agents, equilibrium strategy,  naive strategy, precommitted strategy.

\ms
\no{\it 2020 AMS Mathematics subject classification:}  49L20, 91A65

%\vfill\eject

%\tableofcontents

\section{Introduction}
\label{sect-Introduction}
\setcounter{equation}{0} 

For time inconsistent optimal control problems, different types of agents would choose different strategies, for example the precommitted strategy, the naive strategy, and the equilibrium strategy, see e.g. the survey paper Strotz \cite{Strotz}. Among them, the equilibrium approach by sophisticated agents has received very strong attention. In particular, in the past decade this approach has been extended to continuous time models by many authors, see e.g. 
Bj\"{o}rk-Khapko-Murgoci \cite{BKM}, Bj\"{o}rk-Murgoci-Zhou \cite{BMZ}, Ekeland-Lazrak \cite{EL},  and Yong \cite{Yong}, to mention a few. We also refer to the recent paper Hern\'{a}ndez-Possama\"{i} \cite{HP} for a nice literature review. In this approach, the sophisticated agent will play a game with (infinitely many) future selves, and the goal is to find an equilibrium which is suboptimal and is time consistent in certain sense. 

In this short note, we construct an example in the deterministic continuous time framework such that the unique equilibrium is not Pareto optimal. To be precise, our optimal control problem has a unique equilibrium $\a^*$, but we can construct another control $\hat \a$ such that
\bea
\label{Pareto}
J(t, \hat\a) <  J(t, \a^*) \q\mbox{for all}~t<T,\q\mbox{and}\q J(T, \hat\a) = J(T, \a^*),
\eea
where $J(t,\a)$ is  the dynamic cost function with control $\a$. This raises the serious question on the rationale of using the equilibrium $\a^*$. 

In a non-cooperative game, it is not surprising that an equilibrium may not be Pareto optimal, for example in the well known prisoner's dilemma, cf. Nash \cite[Example 2]{Nash} and Lacey \cite{Lacey}. In that case, since the players do not play cooperatively,  typically due to lack of mutual trust, they may still choose the equilibrium. For the sophisticated agent in our time inconsistent problem, however, the agent is ``playing" the game with future selves, there is no reason that the agent would ``play" in a non-cooperative way. In particular, in our example, it is not rational or say not wise for the agent to choose $\a^*$ instead of $\hat \a$. We shall point out though, that we are not claiming any optimality of the constructed $\hat \a$.
% is optimal in any sense. % and thus the desired control the agent should use. 

It will be very interesting to explore possible alternative notions of equilibrium, or  of ``good" strategies, which we shall leave to future research. To our opinion a good strategy in a dynamic approach should satisfy at least two basic properties: (i) time consistency; and (ii) Pareto optimality. The time consistency has been a natural consideration for time inconsistent problems, as \cite[pp. 173]{Strotz} points out that, when there is intertemporal conflict (namely time inconsistency), the player's problem  ``{\it is then to find the best plan among those that he will actually follow}". The Pareto optimality is a basic requirement in copperative game theory, in fact it is exactly the idea of the notion core, cf. Gillies \cite{Gillies}. In this short note we want to bring into attention the Pareto optimality, which seems less addressed for time inconsistent problems. We note that the precommitted strategy is by definition Pareto optimal but time inconsistent (unless the original problem is time consistent);  and the equilibrium strategy is time consistent but may not be Pareto optimal, as showed in this note. We also note that the naive strategy, while less interesting and less popular in the literature, is also time consistent, see e.g. the recent paper Chen-Zhou \cite{CZ}. However, in the same spirit of Remark \ref{rem-discrete}, we can easily construct an example such that the naive strategy is not Pareto optimal. 
 
We remark that the time consistency of strategies relies on the criterion we take, for example  the equilibrium strategy and the naive strategy satisfy time consistency in different senses. When exploring alternative good strategies, it will be a crucial and intrinsic component to specify this criterion, which in practice relies on the agent's preference.  We would also like to mention the dynamic utility function in \cite{KMZ} and the moving scalarization in Feinstein-Rudloff \cite{FR}, see also the recent paper \.I\c{s}eri-Zhang \cite{IZ}, where the precommitted strategy satisfies both the time consistency and the Pareto optimality. In this approach the utility function for the subproblem over $[t, T]$ is modified and thus it leads to a different function $J(t, \a)$ when $t>0$. Similarly, whether or not to use the dynamic utility function or the original utility function relies on the agent's preference. 

The rest of the paper is organized as follows. In Section \ref{sec:counterexample} we construct the example and prove \reff{Pareto}. In Section \ref{sect-inconsistent} we verify that our example is indeed time inconsistent. Finally in Section \ref{sect-naive} we construct another example where the naive strategy is not Pareto optimal.

\section{An example} 
\label{sec:counterexample}
\setcounter{equation}{0}
Set the time horizon $[0, T]$ with $T=1$, and denote 
\bea
\label{stn}
t_n := 1- 2^{-n}\q \mbox{and}\q s_n := {1\over 4}[t_n + 3 t_{n+1}],\qq n \ge 0.
\eea
 The admissible control set $\cA$ consists of Borel measurable functions $\a: [0, 1]\to [-1, 1]$. 
Consider the following deterministic $2$-dimensional backward controlled system:
\bea
\label{Y12a}
\left.\ba{c}
\dis Y^{1,\a}_t := \int_t^T \a_s ds,\q Y^{2,\a}_t := \int_t^T c(s) \a_s Y^{1,\a}_s  ds,\q t\in [0, T],\\
\dis \mbox{where}\q c(t) := \sum_{n=0}^\infty \big[\1_{[t_n, s_n)}(t) + 6 \1_{[s_n, t_{n+1})}(t)\big].
\ea\right.
\eea
Our time inconsistent optimization problem is:
\bea
\label{Vt}
V_t := \inf_{\a\in \cA} J(t,\a),\q\mbox{where}\q J(t,\a):= Y^{2,\a}_t.
\eea
We remark that, besides the hyperbolic discounting, this type of multidimensional optimization problem is also typically time inconsistent, which includes the well known mean variance optimization problem, see Karnam-Ma-Zhang \cite{KMZ} for more discussions. We shall prove rigorously the time inconsistency for this example  in the next section, and here we focus on the equilibrium approach and the Pareto optimality. 

We first recall the notion of equilibrium:
\begin{defn}
\label{defn-equilibrium}
We say $\a^*\in \cA$ is an equilibrium if, for any $t\in [0, T)$ and $\a\in \cA$,
\bea
\label{equilibrium}
\liminf_{\d\downarrow 0} {1\over \d}[J(t, \a\oplus_{t+\d} \a^*) - J(t,\a^*)] \ge 0,\q\mbox{where}\q \a\oplus_{t+\d} \a^*:= \a\1_{[0, t+\d)} + \a^*\1_{[t+\d, T]}.
\eea
\end{defn}

\begin{prop}
\label{prop-equilibrium}
$\a^*\equiv 0$ is the unique equilibrium.
\end{prop}
\proof (i) We first show that $\a^*\equiv 0$ is an equilibrium. Indeed, for any $t\in [0, T)$, by our construction of $c$, there exists $\d_t >0$ such that $c(s) \equiv c(t)>0$ for all $s\in [t, t+\d_t]$. Then, for any $0<\d \le \d_t$ and any $\a\in \cA$, denoting $\a^\d:= \a\otimes_{t+\d} \a^*$,
\bea
\label{J>0}
\left.\ba{c}
\dis Y^{1, \a^\d}_s = \1_{[0, t+\d]}(s) \int_s^{t+\d} \a_r dr,\ms\\ 
\dis J(t, \a^\d) = Y^{2, \a^\d}_t = \int_t^{t+\d} \!\!\!\! c(s) \a_s Y^{1,\a^\d}_s  ds = c(t)\int_t^{t+\d}  \!\!\!\!\a_s \int_s^{t+\d}\!\!\!\! \a_r drds = {c(t)\over 2} \big(\int_t^{t+\d}\!\!\!\!\a_rdr\big)^2 \ge 0.
\ea\right.
\eea
It is obvious that $J(t, \a^*) = 0 \le J(t, \a^\d)$, then \reff{equilibrium} holds true. In fact, $\a^*$ is an equilibrium in a stronger sense as in He-Jiang \cite{HJ} and Huang-Zhou \cite{HZ}.

(ii) We next show the uniqueness.  Let $\a^*\in \cA$ be an arbitrary equilibrium. For any $t\in [0, T)$ and $\a\in \cA$, let $\d_t>0$ and $\a^\d$ be as in (i). Then for any $0<\d\le \d_t$, noting that $|\a|, |\a^*|\le 1$,
\bea
\label{J-J}
J(t,\a^\d) - J(t, \a^*) = c(t) \int_t^{t+\d} \!\!\big[\a_s Y^{1,\a^\d}_s - \a^*_s Y^{1,\a^*}_s\big]ds =  c(t) Y^{1, \a^*}_t \int_t^{t+\d} \!\! [\a_s  - \a^*_s ]ds + O(\d^2).
\eea
Thus, by setting $\a \equiv -1$ when $Y^{1, \a^*}_t >0$ and $\a \equiv 1$ when $Y^{1, \a^*}_t < 0$, it follows from \reff{equilibrium}  that
\bea
\label{a*1}
\left.\ba{c}
\dis \liminf_{\d\downarrow 0}{1\over \d}  \int_t^{t+\d} [-1  - \a^*_s ]ds \ge 0,\q\mbox{whenever} ~Y^{1, \a^*}_t >0,\\
\dis \limsup_{\d\downarrow 0}{1\over \d}  \int_t^{t+\d} [1  - \a^*_s ]ds \le 0,\q\mbox{whenever} ~ Y^{1, \a^*}_t <0.
\ea\right.
\eea
We recall again that $|\a^*|\le 1$, and that $Y^{1, \a^*}_t$ is continuous in $t$. Then clearly
\bea
\label{a*2}
\left.\ba{c}
\dis \a^*_t = -1\q\mbox{for Leb-a.e.}~ t\in \{s\in [0, T]: Y^{1, \a^*}_s >0\},\\
\dis \a^*_t = 1\q\mbox{for Leb-a.e.}~ t\in \{s\in [0, T]: Y^{1, \a^*}_s < 0\}.
\ea\right.
\eea
Note further that 
\bea
\label{a*3}
\int_0^T \a^*_t Y^{1,\a^*}_t dt = {1\over 2} \Big(\int_0^T \a^*_tdt\Big)^2 \ge 0.
\eea
 This, together with \reff{a*2}, implies that $ Y^{1, \a^*}_t \equiv 0$ for all $t\in [0, 1]$. Therefore, $\a^* \equiv 0$.
 \qed

\begin{prop}
\label{prop-eg}
The following $\hat \a$ dominates $\a^*\equiv 0$ in the sense of \reff{Pareto}:
\bea
\label{hata}
\hat \a_t := \sum_{n=0}^\infty \Big[ \1_{[t_n, s_n)}(t) - \1_{[s_n, t_{n+1})}(t)\Big].
\eea 
\end{prop}
\proof First, we note that
\beaa
s_n-t_n = {3\over 2^{n+3}},\q t_{n+1}-s_n = {1\over 2^{n+3}}.
\eeaa
For any $n$ and $t\in [s_n, t_{n+1})$, we have
\bea
\label{hatY1}
Y^{1,\hat\a}_t &=& -(t_{n+1}-t) + \sum_{m=n+1}^\infty \big[(s_m-t_m) - (t_{m+1}-s_m)\big] =  -(t_{n+1}-t) + \sum_{m=n+1}^\infty {1\over 2^{m+2}}\nonumber\\
 &=& {1\over 2^{n+2}} - (t_{n+1}-t) \ge {1\over 2^{n+2}} - (t_{n+1}-s_n) = {1\over 2^{n+3}} >0; 
\eea
and for $t\in [t_n, s_n)$,
\bea
\label{hatY2}
Y^{1,\hat\a}_t = Y^{1,\hat\a}_{s_n} + (s_n-t) =  {1\over 2^{n+3}} + (s_n-t) >0.
\eea
Then, for any $n$, recalling the $c(t)$ in \reff{Y12a} and $\hat \a$ in \reff{hata},
\bea
\label{hatJ1}
Y^{2,\hat\a}_{t_n} &=& \sum_{m=n}^\infty \Big[\int_{t_m}^{s_m} Y^{1,\hat\a}_t dt - 6 \int_{s_m}^{t_{m+1}} Y^{1,\hat\a}_t dt\Big] \nonumber\\
 &=&  \sum_{m=n}^\infty \Big[\int_{t_m}^{s_m} \big[{1\over 2^{m+3}} + (s_m-t)\big] dt - 6 \int_{s_m}^{t_{m+1}} \big[ {1\over 2^{m+2}} - (t_{m+1}-t)\big] dt \Big] \nonumber\\
 %&=&  \sum_{m=n}^\infty \Big[ \big[{1\over 2^{m+3}}(s_m-t_m) + {1\over 2}(s_m-t_m)^2\big]  - 6 \big[ {1\over 2^{m+2}}(t_{m+1}-s_m) - {1\over 2}(t_{m+1}-s_m)^2\big] \Big] \nonumber\\
&=&  \sum_{m=n}^\infty \Big[ \big[{1\over 2^{m+3}}\times {3\over 2^{m+3}}  + {1\over 2}({3\over 2^{m+3}})^2\big]  - 6 \big[ {1\over 2^{m+2}}\times {1\over 2^{m+3}} - {1\over 2}({1\over 2^{m+3}})^2\big] \Big] \nonumber\\
&=& -\sum_{m=n}^\infty {3\over 2^{2m+7}} = -{1\over 2^{2n+5}} < 0.
\eea
Recall by \reff{hatY1}, \reff{hatY2} that $Y^{1,\a}_t >0$ for all $t<T$. Then, for each $n$,
\bea
\label{hatJ2}
\left.\ba{lll}
\dis t\in [t_n, s_n): \q Y^{2, \hat\a}_t = Y^{2, \hat\a}_{t_n} - \int_{t_n}^t  Y^{1,\hat\a}_s ds \le Y^{2, \hat\a}_{t_n} <0;\\
\dis t\in [s_n, t_{n+1}):\q Y^{2, \hat\a}_t = Y^{2, \hat\a}_{t_{n+1}} - 6 \int_t^{t_{n+1}}  Y^{1,\hat\a}_s ds \le Y^{2, \hat\a}_{t_{n+1}} <0.
\ea\right.
\eea
Since $J(t, \a^*)=0$, this proves \reff{Pareto} immediately.
\qed

\begin{rem}
\label{rem-discrete}
In a discrete time setting: $0=t_0<\cds<t_n = T$, by definition of equilibrium we must have $J(t_{n-1}, \a^*) \le J(t_{n-1}, \a)$, so there is no $\hat \a$ satisfying \reff{Pareto} at $t_{n-1}$. However, following the same spirit one can easily construct examples such that $J(t_i, \hat \a) < J(t_i, \a^*)$ for all $i=0,\cds, n-2$. Here in continuous time model, the ``last" step vanishes and thus the strict inequality holds for all $t<T$.
\end{rem}

\section{The time inconsistency}
\label{sect-inconsistent}
\setcounter{equation}{0}

We now show that the dynamic optimization problem \reff{Vt} is time inconsistent. 
We first note that  \reff{Vt} admits an optimal control for any fixed $t$.

\begin{prop}
\label{prop-optimal}
For any $t\in [0, T)$, the $V_t$ in \reff{Vt} has an optimal control $\bar \a^*\in \cA$.
\end{prop}

Before we prove this result, we use it to show the time inconsistency.

\begin{prop}
\label{prop-inconsistent}
The  dynamic problem \reff{Vt} is time inconsistent.
\end{prop}
\proof Assume by contradiction that \reff{Vt} is time consistent, then there exists  $\bar\a^*\in \cA$ which is optimal for all $t\in [0, T)$. By Definition \ref{defn-equilibrium} this $\bar \a^*$ is an equilibrium, and thus by Proposition \ref{prop-equilibrium} we must have $\bar \a^* \equiv 0$. However, by Proposition \ref{prop-eg} we see that $\a^*\equiv 0$ is not optimal for all $t<T$, which is the desired contradiction.
\qed

\bs

\no{\bf Proof of Proposition \ref{prop-optimal}.} Without loss of generality, we prove the result only at $t=0$. Let $\cX$ denote the set of functions $X: [0, T]\to \dbR$ such that 
\bea
\label{cX}
|X_t - X_s|\le |t-s|\q\mbox{and}\q X_T = 0,
\eea
and we equip $\cX$ with the uniform norm. Then the set $\cX$ is compact and $\a\in \cA$ has one to one correspondence with $X\in \cX$ in the sense: $X_t = \int_t^T \a_sds$ and $\a_t = -X'_t$. Therefore,
\bea
\label{VX}
\left.\ba{c}
\dis V_0 = \inf_{X\in \cX} \big[ -\int_0^T c(t) X_t X'_t dt\big] = - \sup_{X\in \cX} F_\infty(X),\q\mbox{where}\\
\dis F_n(x) :=  \int_0^{t_n} c(t) X_t X'_t dt,\q F_\infty(X) :=  \int_0^T c(t) X_t X'_t dt.
\ea\right.
\eea
Note that
\bea
\label{Fn}
F_n(X) &=& \sum_{m=0}^{n-1} \Big[ c(t_m) \int_{t_m}^{s_m} X_t X'_t dt + c(s_m) \int_{s_m}^{t_{m+1}} X_t X'_t dt\Big]\nonumber\\
&=&\sum_{m=0}^{n-1} \Big[ {c(t_m)\over 2}(X_{s_m}^2 - X_{t_m}^2) + {c(s_m)\over 2}(X_{t_{m+1}}^2-X_{s_m}^2)\Big].
\eea
It is obvious that $F_n$ is continuous in $X$. Moreover,
\bea
\label{Fninfty}
\left.\ba{c}
\dis \sup_{X\in \cX} \big|F_\infty(X) - F_n(X)\big| \le \int_{t_n}^T c(t) |X'_t|\int_t^T |X'_s|ds dt  \\
\dis \le \int_{t_n}^T c(t) (T-t) dt \le {C\over 2^{2n}}\to 0,\q \mbox{as}~n\to\infty.
\ea\right.
\eea
Then $F_\infty$ is also continuous, and thus, by the compactness of $\cX$ there exists $X^*\in \cX$ such that $F_\infty(X^*) = \sup_{X\in \cX} F_\infty(X)$. Therefore, $\bar\a^*_t := -{d\over dt} X^*_t$ is an optimal control for $V_0$.
\qed

\section{Non-optimality of naive strategy}
\label{sect-naive}
In this section we investigate briefly the naive strategy, which is much less popular in the literature. It is well understood that a naive strategy may not be optimal, except in the last step in a discrete time model, as in Remark \ref{rem-discrete}. We now provide an example in continuous time framework such that the naive strategy is Pareto dominated by another strategy. 

Let $\cA$ consist of Borel measurable functions $\a: [0, T]\to \dbR$. 
Consider
\bea
\label{Vt-naive}
V_t := \inf_{\a\in \cA} J(t,\a),\q\mbox{where}\q J(t,\a):= \int_t^T |\a_s - K(t,s)|ds,\q K(t,s) := 2(s-t).
\eea
It is obvious that, for each $t\in [0, T]$, the optimization problem \reff{Vt-naive} on $[t, T]$ has a unique optimal control $\a^t_s := K(t,s)$, $s\in [t, T]$. In particular, for $t_1 < t_ 2$,
\bea
\label{naive-inconsistent}
\a^{t_1}_s = K(t_1, s) \neq K(t_2, s) = \a^{t_2}_s,\q s\ge t_2.
\eea
That is, the dynamic problem \reff{Vt-naive} is time inconsistent. The naive strategy is defined as:
\bea
\label{a*-naive}
\a^*_t := \a^t_t = K(t,t)=0,\q 0\le t\le T,
\eea
and thus
\bea
\label{Ja*-naive}
 J(t,\a^*) =  \int_t^T |\a^*_s - K(t,s)|ds = \int_t^T 2(s-t) ds = (T-t)^2.
\eea

We now set 
\bea
\label{hata-naive}
\hat \a_t := T-t, \q 0\le t\le T.
\eea
Then, for any $0\le t\le T$,
\bea
\label{Jhata-naive}
J(t, \hat\a) =  \int_t^T |\hat\a_s - K(t,s)|ds =  \int_t^T |T+2t - 3s|ds =\int_0^{T-t} |T-t - 3s|ds = {5\over 6}(T-t)^2.
\eea
Thus $\hat \a$ Pareto dominates the naive strategy $\a^*$ in the sense of \reff{Pareto}.

\end{document}